\documentclass[10pt]{amsart}
\usepackage{amsmath,amsthm,latexsym,amsfonts,amssymb,amsthm}
\usepackage{verbatim}
\usepackage{hyperref}
\usepackage{url}

\usepackage{color}
\usepackage{graphicx}

\def\sub#1{_{_{#1}}}
\def\sup#1{^{^{#1}}}

\def\mf#1{\mathfrak{#1}}
\def\mb#1{\mathbb{#1}}
\def\mc#1{\mathcal{#1}}

\def\fmod#1{f\sup{\mc{#1}}}
\def\Rmod#1{R\sup{\mc{#1}}}

\def\mod#1#2{{#1}\sup{\mc{#2}}}

\def\forza#1{\Vdash\sub{#1}}

\def\coLim#1{\underset{#1}{\text{\rm coLim}}\ }

\def\R#1{\mb{R}\sup{#1}}                      
\def\N#1{\mb{N}\sup{#1}}

\def\Im{\text{\rm Im}}

\def\Cond{\Rightarrow}
\def\Recip{\Leftarrow}
\def\Iff{\Leftrightarrow}

\def\ARROW#1{
\text{\begin{picture}(35,18)(15,15)         
            \put(30,23){$_{#1}$}
            \put(18,18){\vector(1,0){30}}
\end{picture}}}

\newcounter{numero}
\newcommand{\Numero}{\setcounter{numero}{1}(\arabic{numero}) }
\newcommand{\numero}{\addtocounter{numero}{1}(\arabic{numero}) }

\newcounter{letra}
\newcommand{\Letra}{\medskip \setcounter{letra}{1}(\alph{letra}) }
\newcommand{\letra}{\medskip \addtocounter{letra}{1}(\alph{letra}) }

\newcounter{romnumero}

\newcounter{bibnumero}

\newtheorem{teo}{Theorem}[subsection]                  
\newtheorem{lema}[teo]{Lemma}
\newtheorem{prop}[teo]{Proposition}
\newtheorem{cor}[teo]{Corollary}

\def\bteo{\begin{teo}}
\def\eteo{\end{teo}}
\def\bprop{\begin{prop}}
\def\eprop{\end{prop}}
\def\bcor{\begin{cor}}
\def\ecor{\end{cor}}
\def\blema{\begin{lema}}
\def\elema{\end{lema}}

\theoremstyle{remark}                               
\newtheorem{obs}[teo]{Remark}

\def\bobs{\begin{obs} }
\def\eobs{\end{obs}}
\newenvironment{dem}{ [{\it Proof\/}]\rm\hskip3mm }{\hfill$\square$\vskip5mm}       
\def\bdem{\begin{dem}}
\def\edem{\end{dem}}

\begin{document}

\author{Gabriel Padilla}
\author{Andr\'es Villaveces}
\email{gipadillal@unal.edu.co, avillavecesn@unal.edu.co}

\address{Departamento de Matem\'aticas\\ Universidad Nacional de Colombia.\\AK30 Cl45 Edif.404\\Bogot\'a 111321. Colombia\\ Tel. +5713165207}

\title[Non standard Cohomology and Generic Models]{Non standard Cohomology for Equivariant Sheaves: the role of Generic Models}



\begin{abstract}
		We introduce a new class of generic cohomologies and 
		show how, in some cases, they simplify non standard
        cohomologies \cite{brunjes,macintyre2}.  
        For doing so, we use a previous generalization of the Generic Model Theorem 
		for equivariant exact  presheaves of structures; extending
		the results of Macintyre and Caicedo \cite{caicedo,macintyre}.

\end{abstract}
\maketitle

\section*{Foreword}\label{section foreword}
We work with an equivariant version of the sheaves of structures
introduced by Comer~\cite{comer} and Macintyre~\cite{macintyre} and
later further expanded by Caicedo~\cite{caicedo}, suitable for
working with a transformation group $G$ acting on the sheaf.
A $G$-structure is a structure with an action such that 
elements of $G$ commute with language symbols.
Given a presheaf of $G$-structures $\mc{M}$ on $X$; through a previous extension of (topo)logical
``truth'' due to Caicedo \cite{caicedo} to the equivariant context,  we show that
\begin{itemize}
 \item If $\mc{M}$ is exact, then it has a generic $G$-model.
 \item If $(\mc{M},d)$ is a differential presheaf of structures with generic model
 $\mc{M}\sup{gen}$ then there is a canonical generic cohomology structure
 $H\sup{gen}(X,\mc{M})$ with nice category-theoretic properties.
\end{itemize}
This is done through a
simplification in the presentation of the pointwise forcing
relation. The article follows this sequence: Sections 1 and 2 provide the basics on $G$-structures and  
equivariant presheaves with fibers of that kind. In section 3
we review the local semantics and their behavior. 
Section 4 is devoted to the construction of equivariant generic models  
and what we mean by their {\it ``genericity''}
(Theorem \ref{teo genericidad}), so we extend both the Generic Model Theorem 5.2 of
\cite{caicedo} and Theorem 3 of \cite{macintyre}
to the equivariant context. 
We also introduce a new class of cohomologies related to differential
generic models.
Our aim here is therefore to establish the first results towards
a Model Theoretic analysis of geometric structures beyond
sheaves.\\[2mm] 
From now on, we fix a group $G$.

\section{$G$-structures}\label{sssection traduccion por submersion}
We follow usual model theory conventions (see for example
\cite{hodges,marker}) in our definition (and  notation) of a {\bf language} 
$\mc{L}=(\mc{F},\mc{R},\mc{C})$ and a 
{\bf structure} $\mc{M}=\left(M,\mc{F}\sup{\mc{M}},\mc{R}\sup{\mc{M}},\mc{C}\sup{\mc{M}}\right)$.
Relations, functions, constants, arities, formulas and\linebreak  semantics (such as 
$\mc{M}\models\varphi(a)$) are understood on this context.
\subsection{}\label{sssection morphisms of structures}
A {\bf morphism} of structures $\mc{M}\ARROW{\alpha}\mc{N}$
is a map between universe sets $M\ARROW{\alpha}N$  that commutes with the language symbols. 
In particular, we only ask
$\alpha\left(R\sup{\mc{M}}\right)\subset R\sup{\mc{N}}$ for each $R\in\mc{R}$. 
Geometric reasons for this will be appreciated soon. 
An {\bf isomorphism} is a bijective morphism whose inverse is also a morphism. 

\subsection{}\label{sssection transfilled morphisms}
A morphism $\mc{M}\ARROW{\alpha}\mc{N}$ is {\bf transfilled} (short
for {\it ``transversally filled''}) iff
$\alpha\sup{-1}\left(R\sup{\mc{N}}\right)\subset R\sup{\mc{M}}$ for
each relation symbol $R\in\mc{R}$. An {\bf embedding} (resp. a {\bf
  submersion}) is an injective (resp. surjective) transfilled
morphism.
Given two structures $\mc{M,N}$ such that $M\subset N$; we say that
$\mc{M}$ is a {\bf substructure} of $\mc{N}$ iff the inclusion map is
an embedding, in
that case we write $\mc{M}\leq\mc{N}$.\\[2mm]
Given a transfilled morphism $\mc{M}\ARROW{\alpha}\mc{N}$; there is a unique substructure
$\mc{I}(\alpha)$ of $\mc{N}$ whose universe $\Im(\alpha)$ is the image set of  $\alpha$. 
From the obvious equivalence relation on $M$ induced by $\alpha$ we also obtain a quotient model $\mc{M}/\sim$ whose universe is the quotient set $M/\sim$; 
the quotien projection $\mc{M}\ARROW{q}\mc{M}/\sim$ is a submersion and the induced arrow $\mc{M}/\sim\ARROW{\overline{\alpha}}\mc{I}(\alpha)$ 
is an isomorphism.

\subsection{}\label{sssection G-estructuras}
Fix some group $G$. A {\bf $G$-structure} is a structure $\mc{M}$ such that
\begin{enumerate}
	\item The universe $M$ is endowed with an action $G\times M\ARROW{\Phi}M$.
	\item The action commutes with the language symbols. More precisely
	\begin{enumerate}
	 \item \underline{The set of constants is invariant}: $G\mod{\mc{C}}{A}=\mod{\mc{C}}{A}$.
	 \item \underline{Relations are invariant subsets}: $G\sup{n_R}\Rmod{A}=\Rmod{A}$ 
	 for each $R\in\mc{R}$. In other words; given $\left(x\sub1,\dots,x\sub{n_R}\right)\in\Rmod{A}$
	 and $g\sub1,\dots,g\sub{n_R}\in G$, we also have
	 $\left(g\sub1 x\sub1,\dots,g\sub{n_R} x\sub{n_R}\right)\in\Rmod{A}$. 
	 \item \underline{ functions are $G$-equivariant}: For each $f\in\mc{F}$ with arity $n$,
	$g\in G$ and $x\sub1,\dots,x\sub{n}\in A$;
	\[
      \fmod{A}(g x\sub1,\dots,g x\sub{n})=
      g\fmod{A}(x\sub1,\dots,x\sub{n})
	\]
	 
	 	\end{enumerate}
\end{enumerate}

\subsection{Examples}\label{ssection examples} Here there are some examples: \\[2mm]
\Numero A module $\mc{M}=\left(M,+,0\sub{\mc{M}}\right)$ over a principal ideal domain $(D,+,\cdot,0,1)$ is a $(D,+,0)$-structure. The equivariance of the map $+$ implies linearity of scalar stretching since 
$am+am'=+(am,am')=a\cdot[+(m,m')]=a(m+m')$ for $a\in D$, $m,m'\in M$.  Also 
$D\cdot0\sub{\mc{M}}=0\sub{\mc{M}}$ so the set of constants is invariant.\\[2mm]
\numero For a compact group $G$ and a topologic $G$-space $X$; each $g\in G$ provides, by left multiplication, a homeomorphism $X\ARROW{g}X$. These induce a family of chain isomorphism on singular chains
\[SC\sub{*}(X)\ARROW{g}SC\sub{*}(X).\]
Notice that $G$ acts linearly on the singular chain groups $SC\sub{*}(X)$ and also on the homology groups
$H\sub{*}(X)$ so these cases are similar to example (1). The same can be done for Lie groups, smooth manifolds,
smooth forms and De Rham cohomology. It actually can be extended to more complicated 
(co)homology theories, see for instance \cite{illinois}.\\[2mm]
\numero  This notation will be useful in the sequel: For any set $Y$ and a cardinal $\kappa\leq|Y|$, write
$Y\sup{[\kappa]}=\{A\subset Y:|A|=\kappa\}$. The sets $Y\sup{[<\kappa]}$
and $Y\sup{[\leq\kappa]}$ are defined in a similar way.
A countable polyhedron is a pair $(\mc{S},\subset)$ in the language of posets 
$\mc{L}=(<)$, where $\mc{S}\subset \mb{N}\sup{[<\infty]}$ is a family of finite subsets of the set of natural numbers
$\mb{N}$ in  which is hereditary in the following sense:
if $u\subset v\in\mc{S}$ then $u\in\mc{S}$; see \cite{mijares} for details.
The standard $n$-simplex corresponds
to $\Delta\sup{n}=\left(n,\mc{P}(n)\right)$ where, as usual, $0=\emptyset$
$n=\{0,\dots,n-1\}$ for $n>0$, and $\mc{P}(n)$ is the set parts of $n$. For this finite polyhedron, the unique
subgroup $G\leq S\sub{n}$ such that $\Delta\sup{n}$ is a $G$-structure is the trivial group $G=\{e\}$. 
On the other hand; the geometric boundary $\partial\left(\Delta\sup{n}\right)=\mc{P}(n)\backslash\{n\}$, is a 
$G$-structure for any subgroup of $G\leq S\sub{n}$.\\[2mm]
\numero We due the following example/remark to X. Caicedo. 
The orbit set $M/G$ of a $G$-structure $\mc{M}$ is not necessarily a structure in the same 
language. For this to happen, we should consider a stronger notion $G$-structure,
replacing 2(c) by  more restrictive condition as, for instance, that each function $\fmod{M}$ to be coordinatewise equivariant; then the induced quotient functions
would make sense in $\left(\mc{M}/G\right)\sup{n}$. It is possible to construct generic orbit models coming from sheaves
of {\it "strong"} $G$-structures; however the first two examples of this list would be excluded.

\subsection{}\label{ssection orbits and stabilizers}
A {\bf morphism} of $G$-structures is a $G$-equivariant morphism. 
A {\bf $G$-substructure} of $\mc{M}$ is a $G$-invariant substructure. 
The composition of $G$-equivariant morphisms 
(resp. embeddings, submersions, etc.) is an
arrow of the same kind.

\subsection{}\label{ssection notational convention}\label{ssection limites de estructuras}
The family of $G$-structures and $G$-equivariant morphisms
(resp. embeddings, submersions, elementary embeddings) 
is a category, we will denote it by $\mf{M}\sub{G}$ (resp. $\mf{M}\sup{\leq}\sub{G}$, 
$\mf{M}\sup{\geq}\sub{G}$, $\mf{M}\sup{\prec}\sub{G}$).
For an inverse system of $G$-structures
$\{\mc{M}\sub{i}:i\in\mc{D}\}$ and $a\in M\sub{i}$ for some $i\in \mc{D}$ we write $[a]$ for the
germ of $a$ in the colimit
$\mc{M}=\coLim{i\in\mc{D}}\mc{M}\sub{i}$. The limit action of $G$ on
$\mc{M}$ is well defined.

\bprop\label{cor limites de estructuras}
	Let $\mc{M}=\coLim{i\in\mc{D}}\mc{M}\sub{i}$ be a colimit in $\mf{M}\sub{G}\sup{\leq}$, 
	and $a\in M\sup{n}\sub{i}$ for some $i\in\mc{D}$. If 
	$\varphi(v)$ has no $\neg,\forall$; then
	$\mc{M}\models\varphi([a])$ if and only if there is some $j\leq i$ such that 
	$\mc{M}\sub{j}\models\varphi(\rho\sub{ji}(a))$.
\eprop
\bdem $(\Cond)$ By induction on formulas. For instance:
	\begin{itemize}
			\item \underline{$\varphi$ is $t(v)=s(v)$:}
			If $\mc{M}\models\varphi([a])$ then 
			$t\sup{\mc{M}}([a])=s\sup{\mc{M}}([a])$
			so $\left[t\sup{\mc{M}\sub{i}}(a)\right]=\left[s\sup{\mc{M}\sub{i}}(a)\right]$.
			There is some $j\leq i$ such that 
			$\rho\sub{ji}\left(t\sup{\mc{M}\sub{i}}(a)\right)=
			\rho\sub{ji}\left(s\sup{\mc{M}\sub{i}}(a)\right)$. Since $\rho\sub{ji}$
			is a morphism it commutes with terms, so\linebreak    
			$t\sup{\mc{M}\sub{j}}(\rho\sub{ji}(a))=s\sup{\mc{M}\sub{j}}(\rho\sub{ji}(a))$;
			therefore $\mc{M}\sub{j}\models\varphi(\rho\sub{ji}(a))$.		
	\end{itemize} 
	The case $\varphi(v):=[t(v)\in R]$ is similar.
	Next suppose that
	\begin{itemize}
			\item \underline{$\varphi(v)$ is
                          $\psi(v)\wedge\theta(v)$:}
			By induction, assume the statement for both
                        $\psi(v)$ and $\theta(v)$. 
			$\mc{M}\models\varphi([a])$ iff
                        $\mc{M}\models\psi([a])$ and
                        $\mc{M}\models\theta([a])$.
			Then $\exists k,k'\leq i$ such that 
			$\mc{M}\sub{k}\models\psi(\rho\sub{ki}(a))$ and 
			$\mc{M}\sub{k'}\models\theta(\rho\sub{k'i}(a))$.
			Take any $j\leq k,k'$ and 
			$\rho\sub{ji}(a)=\rho\sub{k'j}(\rho\sub{k'i}(a))=
			\rho\sub{kj}(\rho\sub{ki}(a))$. Now, notice
                        that all restrictions maps are embeddings
			and, by our assumptions, we can
                        suppose that $\psi,\theta$ are quantifier-free
                        formulas.
			Since $\rho\sub{kj}$ is an embedding and
                        $\mc{M}\sub{k}\models\psi(\rho\sub{ki}(a))$,
                        we deduce that
			$\mc{M}\sub{j}\models\psi(\rho\sub{ji}(a))$. Similarly,
			we get
                        $\mc{M}\sub{j}\models\theta(\rho\sub{ji}(a))$.\linebreak  Therefore,
			$\mc{M}\sub{j}\models\varphi(\rho\sub{ji}(a))$.
	\end{itemize} 
	The other cases are similar. The converse $(\Recip)$ holds
        because morphisms preserve the validity of
	formulas without $\neg,\forall$.
\edem

\section{Sheaves of $G$-structures}\label{section haces de estructuras}

\subsection{}\label{ssection haces de estructuras}\label{prehaces de estructuras}
The definitions (and notation) of presheaves, restrictions, and
induced sheaves are directly taken from \cite{bredon,godement}.
A {\bf presheaf of $G$-structures} on $X$ is a presheaf
$(\mc{T},\subset)\ARROW{\mc{M}}\mf{M}\sub{G}$. 
Each open subset $U$ of  $X$ is sent to some $G$-structure
$\mc{M}\sub{U}$ and each inclusion of open subsets $U\subset V$ is
mapped to the corresponding\linebreak  equivariant restriction morphism
$\mc{M}\sub{V}\ARROW{\hskip-2mm \rho\sub{UV}}\mc{M}\sub{U}$. When $G$
is trivial we talk about a ``presheaf of structures''. 

\subsection{Examples}\label{ssection ejemplos prehaces de G-estructuras} 
The following examples were inspired by Gendron~\cite{gendron}; we
further develop them later in this paper as examples~
\ref{ssection ejemplo modelo generico del haz RS} and~\ref{ssection
  non standard cohomology}.
Let's consider some examples on $X=\N{}$ the set of natural numbers with the discrete topology.\\[2mm] 
  \Numero The sheaf of real sequences $\mc{RS}$ is given as follows: For    
  each $U\subset\N{}$ define $\mc{RS}\sub{U}=\R{U}$ as the set of all maps from $U$ to $\R{}$.
  For each inclusion $U\subset V$ there is a restriction map 
  $\mc{RS}\sub{V}\ARROW{\hskip-2mm\rho_{UV}}\mc{RS}\sub{U}$ given by 
 $\alpha\mapsto\alpha|\sub{U}$.
  Coherence and exactness are straightforward. The group
  $\mb{Z}$ of integer numbers acts on $\mc{RS}$ with the translations induced by its structure as
  an additive subgroup of $\R{}$, more precisely
  $\mb{Z}\times\mc{RS}\sub{U}\ARROW{}\mc{RS}\sub{U}$
  is given by $(n,\alpha)(i)=n+\alpha(i)$, $\forall i\in U$.\\[2mm]
 \numero The presheaf $\mc{G}$ of graphs on $\N{}$
 is given as follows: For each $U\subset\N{}$ define 
$\mc{G}\sub{U}=2\sup{U^{[2]}}$.
 An element of $\mc{G}\sub{U}$ is a function
 $U\sup{[2]}\ARROW{\alpha}2$ that decides, for each possible edge
 $e=\{u,v\}\subset U$, whether $u,v$ are connected
 ($\alpha(\{u,v\})=1$) or not ($\alpha(\{u,v\})=0$). 
 If $U\subset V$ then $U\sup{[2]}\subset V\sup{[2]}$, so the restriction
   $\rho\sub{UV}(\alpha)=\alpha|\sub{U^{[2]}}$
  makes sense: it corresponds 
 to the graph obtained by dropping the vertices in $U\backslash V$  \cite{diestel}. 
 This presheaf is exact but not coherent.  It is possible to define in a similar way 
 a presheaf $\mc{P}\sup{k}$ of $k$-polyhedra, where $k>0$ is the
 geometric dimension allowed; as
  $k$ grows, there are more possibilities
  to extend local $k$-polyhedra, so $\mc{P}\sup{k}$ 
  becomes less coherent. 
  
\section{Local Semantics}\label{section local semantics}
\subsection{Pointwise semantics}\label{ssection point forcing}
Fix some presheaf of $G$-structures $\mc{M}$ on $X$ and a point $x\in X$. 
Let $\varphi(v)$ be a formula in free variables
$v=(v\sub1,\dots,v\sub{n})$. We also fix (temporarily) 
an open set $U\ni x$
and some element $a\in\mc{M}\sub{U}$; we say that $\mc{M}$ {\bf forces}
$\varphi(a)$ at $x$, and write $\mc{M}\forza{x}\varphi(a)$,
according to the following induction:\\[2mm]
\Numero \underline{$\varphi(v)$ has no $\neg,\forall$}: Then
$\mc{M}\forza{x}\varphi(a)$ iff there is some open set $V\subset V\ni$
with $x\in V$, 
such that $\mc{M}\sub{V}\models\varphi\left(\rho\sub{UV}(a)\right)$. Notice that, by Proposition \ref{cor limites de estructuras}
this is equivalent to require that $\mc{M}\sub{x}\models\varphi\left([a]\sub{x}\right)$.\\
\numero \underline{$\varphi(v)$ is $\neg\psi(v)$}:
$\mc{M}\forza{x}\varphi(a)$ iff there is some
open set $V\subset U$ with $V\ni x$ such that $\mc{M}\not\Vdash_{y}\psi(a)$ for
all $y\in V$.\\
\numero \underline{$\varphi(v)$ is $\psi(v)\rightarrow\nu(v)$}:
$\mc{M}\forza{x}\varphi(a)$ iff there is some open set $V\subset U$
with $V\ni
x$ such that, for all $y\in V$; if
		$\mc{M}\forza{y}\psi(a)$ then $\mc{M}\forza{y}\nu(a)$.\\
\numero \underline{$\varphi(v)$ is $\forall w\psi(v,w)$}:
$\mc{M}\forza{x}\varphi(a)$ iff  there is some open set $V\subset U$
with $\ni x$ such that,
for each $y\in V$ and each $b\in\mc{M}\sub{V}$, we have
$\mc{M}\forza{y}\psi\left(a,b\right)$.

\bprop\label{lema forcing for presheaves}
      Pointwise semantics is equivalent to Caicedo's local\linebreak semantics at \cite{caicedo}.
\eprop
\bdem
	See Proposition 4.1.1; \cite[p.6]{padilla2}.
\edem

\subsubsection{}\label{ssection properties point forcing}
Given a point $x\in X$, an open set $U\ni x$, a formula $\varphi(v)$ in free variables 
$v=(v\sub1,\dots,v\sub{n})$ and some $a\in M\sub{U}\sup{n}$; the following properties hold:  
\begin{itemize}
	\item[\Letra] \underline{Local Semantics}: $\mc{M}\forza{x}\varphi(a)$ iff there is some open set 
	$U\supset V\ni x$  such that $\mc{M}\forza{y}\varphi(a)$ for all $y\in V$.
	\item[\letra] \underline{Classical Semantics}: For  an
          isolated point $x\in X$ we have that
	$\mc{M}\forza{x}\varphi(a) \Iff \mc{M}\sub{x}\models\varphi([a]\sub{x})$.			
	\item[\letra] \underline{Excluded Middle Principle}: 
	$\mc{M}\forza{x}\forall u\forall v(u=v\vee u\neq v)$
	iff there is an open set $U\ni x$ such that at all elements
        $y\in U$, all pairs of sections $a,b$ defined at $y$ are
        forced at $y$ to be equal or different. When $\mc{M}$ is a
        sheaf and the base space $X$ is Hausdorff, this means exactly
        that the induced sheafspace is Hausdorff \emph{in some
          neighborhood} of $x$.
\end{itemize}

\subsection{Open semantics}\label{sssection local forcing}
Given a presheaf of $G$-structures $\mc{M}$, an open set $U\subset X$ and some 
$a\in\mc{M}\sub{U}$; we say that $\mc{M}$ {\bf forces $\varphi(a)$ in $U$}, and we write
$\mc{M}\forza{U}\varphi(a)$, iff $\mc{M}\forza{x}\varphi(a)$ for all $x\in U$. 
By \S\ref{ssection properties point forcing}-(1),
$\mc{M}\forza{x}\varphi(a)$ $\Iff$ there is some neighborhood $U\supset V\ni x$ such that 
$\mc{M}\forza{V}\varphi(a)$.\\[2mm]
The validity of $\varphi(a)$ is related to the topology of $X$:  
\begin{itemize}
	\item[\Letra] \underline{Restrictions}: If $U\subset V$ then $\mc{M}\forza{V}\varphi(a)\ \Cond\ \mc{M}\forza{U}\varphi(a)$. 
	\item[\letra] \underline{Coverings}: $\mc{M}\forza{U_i}\varphi\left(a|\sub{U_i}\right)\ \forall i\ \Cond 
	\mc{M}\forza{\underset{i}{\cup} U_i}\varphi(a)$.
	\item[\letra] \underline{Existential quantifier}: $\mc{M}\forza{U}\exists\nu\varphi(a,\nu)$ iff there 
	is an open covering\linebreak  $\underset{i}{\cup} U_i\supset U$ and some $b\sub{i}\in\mc{M}\sub{U_i}$ for each
	$i$; such that $\mc{M}\forza{U_i} \varphi\left(a,b\sub{i}\right)$ for each $i$.		
\end{itemize}

\bprop\label{prop maximum principle}	
	{\bf [Maximum principle]} 
	Let $\mc{M}$ be an exact presheaf of $G$-structures on $X$, $x\in X$ a point, 
	$U\ni x$ an open set and and $a$ in $\mc{M}\sub{U}$.
	If $\mc{M}\forza{U}\exists\nu\varphi(a,\nu)$ then there is 
	some open subset $V\subset U$ and $b\in\mc{M}\sub{V}$ such that $U\subset\overline{V}$ and
	$\mc{M}\forza{V}\varphi(a,b)$.
\eprop
\bdem
	This is a translation of Theorem 3.3 in \cite[p.18]{caicedo}.	
	See Proposition  4.2.2; \cite[p.13]{padilla2}. 
\edem

\section{Equivariant generic models}\label{section equivariant generic models}\label{section EGMT}
In this section we show how the models constructed at 
\S\ref{section equivariant generic models} are generic. Fix a presheaf of $G$-structures
$\mc{M}$ on $X$. 

\subsection{}\label{ssection filtros genericos}
For the definitions and existence of filters and ultrafilters  see
\cite[p.83,135]{kelley}. The following definition is from~\cite{caicedo}.
A (non trivial) filter generated by a family of open subsets 
$\mb{F}$ in $X$  is {\bf generic with respect to the presheaf $\mc{M}$ } iff:
\begin{enumerate}
		\item For each formula $\varphi(v)$ in free variables $v=(v\sub1,\dots,v\sub{n})$, 
		$U\in\mb{F}$ and $a\in M\sub{U}\sup{n}$; there is some  
		$U\supset V\in\mb{F}$ such that $\mc{M}\forza{V}\varphi(\sigma)$ or $\mc{M}\forza{V}\neg\varphi(\sigma)$. 
		\item For each formula $\varphi(v,w)$ in free
                  variables $v=(v\sub1,\dots,v\sub{n})$ and  
		$w=(w\sub1,\dots,w\sub{m})$; each $U\in\mb{F}$ and $a\in M\sub{U}\sup{n}$; 
		if $\mc{M}\forza{U}\exists w\varphi(a,w)$ then		
		there is some $U\supset V\in\mb{F}$ and $b\in M\sub{V}\sup{m}$ such that 
		$\mc{M}\forza{V}\varphi(a,b)$.
\end{enumerate}

\bprop\label{prop generic filters2}
      {\bf [Existence of generic filters]} If $\mc{M}$ is an exact presheaf of $G$-structures; then  
      every ultrafilter in $X$ (generated by a family of open subsets) is generic with respect to $\mc{M}$.
\eprop
\bdem
      For condition \S\ref{ssection filtros genericos}-(1) apply the same proof of 
      Theorem 5.1 at \cite[p.27]{caicedo}. For condition
      \S\ref{ssection filtros genericos}-(2) the main argument is the maximum principle which,
      in our context, only requires the hypothesis of exactness.
\edem

\subsection{}\label{ssection modelos genericos}
A {\bf generic model} for a presheaf of $G$-structures $\mc{M}$ is the colimit structure
$\mc{M}\sup{gen}=\coLim{U\in\mb{F}}\mc{M}\sub{U}$ on a generic filter $\mb{F}$. 

\bcor\label{teo existen filtros genericos}\label{teo modelo generico}
    {\bf [Existence of equivariant generic $G$-models]}
    Every exact presheaf of  $G$-structures $\mc{M}$ has a generic $G$-model $\mc{M}\sup{gen}$.
\ecor

\subsection{Example }\label{ssection ejemplo modelo generico del haz RS}
For the sheaf $\mc{RS}$ of real sequences on $\mb{N}$;  
fix an ultrafilter $\mb{F}$ in $\mb{N}$. 
Given $U\subset\N{}$, the quotient projection 
$\mc{RS}\sub{U}\ARROW{\hskip-2mm q_U}\mc{RS}\sup{gen}$
sends each $\alpha$ to its $\mb{F}$-germ $[\alpha]$. 
In particular, the map $\mc{RS}\sub{\N{}}\ARROW{\hskip-2mm q_{\mb{N}}}\mc{RS}\sup{gen}$ is surjective:
 For each $[\alpha]\in\mc{RS}\sup{gen}$, $\alpha\in\mc{RS}\sub{U}$; take 
 $\beta\in\R{\N{}}$ defined as $\beta(i)=\alpha(i)$ if $i\in U$, and $\beta(i)=0$ if $i\not\in U$. 
 Then $[\alpha]=[\beta]=q\sub{\N{}}(\beta)$. We deduce that
 $\mc{RS}\sup{gen}$
  is the quotient of $\R{\mb{N}}$ by the equivalence modulo the ultrafilter
 $\mb{F}$, i.e. the ultraproduct of $\R{}$ which leads to the structure of non standard real numbers 
 $\mc{RS}\sup{gen}=\sup{*}\R{}$. We should also notice that this is a sheaf of strong
 $\mb{Z}$-structures, in the sense of \S\ref{ssection examples}-(4). 
 Then  
 \[
      (\mc{RS}/\mb{Z})\sup{gen}\cong \frac{(\mb{R}/\mb{Z})\sup{\mb{N}}}{\sim\sub{\mb{F}}}=
      \frac{(\mb{S}\sup{1})\sup{\mb{N}}}{\sim\sub{\mb{F}}}=\sup{*}\mb{S}\sup{1}
 \]
 is the non-standard unit circle group. Genericity, in this case, is just the
 universal semantic property of ultraproducts. 

\subsection{}\label{ssection Godel translation}
Let us show the behavior of the forcing relation under double negations.
We start with two easy statements, the proofs are left to the reader
who can go to \cite{caicedo} for more details.

\blema\label{lemma double negation}
    Let $\varphi(v)$ be a positive formula. Then 
	$\mc{M}\forza{U}\neg(\neg\varphi(a))$ iff there is some open set $V\subset U$
	such that $V$ is dense in $U$ and $\mc{M}\forza{V}\varphi(a)$.
\elema

\blema\label{lemma filters and dense subsets}
    Let $\mb{F}$ be a maximal filter of open sets in $X$, and $U\in \mb{F}$.
    If $V\subset U$ is open and dense in $U$, then $V\in \mb{F}$.
\elema

The {\bf G\"odel translation} $\varphi\sub{\mathbb{G}}$ of some formula $\varphi$ is defined, by induction, as follows: 
\begin{itemize}
 \item $\varphi\sub{\mathbb{G}}$ is $\neg(\neg \varphi)$ for an atomic formula $\varphi$.
 \item $(\varphi\wedge\psi)\sub{\mathbb{G}}=\varphi\sub{\mathbb{G}}\wedge\psi\sub{\mathbb{G}}$.
 \item $(\varphi\vee\psi)\sub{\mathbb{G}}=\neg
 \left(\neg\varphi\sub{\mathbb{G}}\wedge\neg\psi\sub{\mathbb{G}}\right)$.
 \item $(\neg\varphi)\sub{\mathbb{G}}=\neg\left(\varphi\sub{\mb{G}}\right)$.
  \item $(\forall v\varphi)\sub{\mb{G}}=\forall v\left(\varphi\sub{\mathbb{G}}\right)$.
 \item $(\exists v\varphi)\sub{\mb{G}}=\neg\forall v\left(\neg\varphi\sub{\mathbb{G}}\right)$. 
\end{itemize}

\bteo\label{teo genericidad}
    {\bf[Equivariant generic model theorem]} Let $\mc{M}$ be a
    sheaf of $G$-structures
    on $X$ and $\mc{M}\sup{gen}$ the generic model induced by some
    generic filter    
    $\mb{F}$ on $X$. For each formula $\varphi(v)$, 
    $U\in\mb{F}$ and $a\in\mc{M}\sub{U}$; the following statements are
    equivalent: 
    \begin{enumerate}
	  \item $\mc{M}\sup{gen}\models\varphi([a])$.
	  \item
            $\mc{M}\forza{V}\varphi\sub{\mb{G}}\left(a\right)$
            for some $U\supset V\in\mb{F}$.
	  \item $\{x\in U:
            \mc{M}\forza{x}\varphi\sub{\mb{G}}(a)\}\in\mb{F}$.
    \end{enumerate}
\eteo
\bdem
 See Theorem 4.3.3 at \cite{padilla2}.
\edem

\section{Generic cohomology}\label{section IGC}

The model theory of cohomology has had some earlier results
(see~\cite{macintyre2}). Our notion is adapted to sheaves and
$G$-sheaves - we exploit the functoriality of generic models in our
definition and provide various examples that show that this\linebreak  extends
classical notions of cohomology, and furthermore allows us to define
new extensions.

\subsection{}\label{ssection generic functoriality}
Let us fix two exact presheaves $\mc{M,N}$ of $G$-structures. Then;
\begin{enumerate}
	\item By colimit properties; $\mc{M}\sup{gen}$ inherits a natural action, so it is a $G$-structure.
	 \item Each morphism of presheaves (i.e. natural transformation) $\mc{M}\ARROW{}\mc{N}$ induces a morphism 
	 $\mc{M}\sup{gen}\ARROW{\tau}\mc{N}\sup{gen}$ of structures.  
	 If the first arrow is equivariant then the second is amorphism of $G$-structures.

\end{enumerate}

\subsection{}\label{ssection generic cohomology}
When $\mf{M}$ an abelian category there is a
zero-structure. A {\bf differential} of $\mc{M}$ is a transfilled endomorphism
$\mc{M}\ARROW{d}\mc{M}$ such that $d\sup2=0$. We define the {\bf generic cohomology
of $X$  with values in $(\mc{M},d)$} as 
\[
	H\sup{gen}(X,\mc{M}):=H\left(\mc{M}\sup{gen},d\right)
	=\frac{\text{ker}\left[\mc{M}\ARROW{d}\mc{M}\right]}{\text{Im}\left[\mc{M}\ARROW{d}\mc{M}\right]}
\]
These {\it ``cohomology structures''} are well defined by
\ref{sssection transfilled morphisms};
they extend the usual   notion of cohomology. Since the space
$X$ belongs to any ultrafilter, there is a   quotient map
$\mc{M}\sub{X}\ARROW{}\mc{M}\sup{gen}$ which commutes with the differential $d$, and
a well defined morphism of cohomology   structures
\[
	H(X,\mc{M})\ARROW{}H\sup{gen}(X,\mc{M})		
\]
 
\subsection{}\label{ssection Example smooth G-manifolds}
Usual examples come from algebraic topology context. Among them,
\begin{enumerate}
	\item Singular homology groups for $G$-spaces \cite{bredon3}.
	\item De Rham cohomology, for smooth $G$-manifolds \cite{bott,bredon,Greub}.
	\item Intersection homology \cite{gmp} and cohomology \cite{illinois}. 
	\item $(N,q)$-cohomologies of amplitude $1\leq k<N$, and $N>0$; \cite{dubois,kapranov}.
\end{enumerate}
The relation between these presheaves is usually given 
given in terms of\linebreak  cohomology spectral sequences.

\subsection{Example: Non standard generic cohomology}\label{ssection non standard cohomology} 
For similar examples  coming from other contexts, see \cite{brunjes}.
Consider the sheaf $\mc{Z}\sub{n}\mc{S}$ of sequences on  (subsets of) $\mb{N}$ with values on 
the ring $\mb{Z}\sub{n}$ of integers modulo some fixed $n>0$;
let $n=p\sub{1}\sup{r_1}\cdots p\sub{s}\sup{r_s}$ be its prime decomposition.\\[2mm]
For each $U\subset\N{}$ define $\mc{Z}\sub{n}\mc{S}\sub{U}=\mb{Z}\sub{n}\sup{U}$ and for each inclusion $U\subset V$ 
  let $\rho_{UV}(\alpha)=\alpha|\sub{U}$ as usual.  
  A linear map $\mb{Z}\sub{n}\sup{\mb{N}}\ARROW{d}\mb{Z}\sub{n}\sup{\mb{N}}$
 will induce a sheaf endomorphism $\mc{Z}\sub{n}\mc{S}\ARROW{d}\mc{Z}\sub{n}\mc{S}$ iff $d$ is diagonal; i. e.
 iff for each $i\in\mb{N}$,  the element ${\tt e}\sub{i}\in \mb{Z}\sub{n}\sup{\mb{N}}$ of the canonical basis satisfies
 $d({\tt e}\sub{i})=\overline{a}\sub{i}{\tt e}\sub{i}$ for some $\overline{a}\sub{i}\in \mb{Z}\sub{n}$. 
 Notice $d\sup{N}(\overline{x})=\left(\overline{a}\sub{i}\sup{N} \overline{x}\sub{i}\right)_{i\in\mb{N}}$		
 for all $x\in \mb{Z}\sup{\mb{N}}$, $N>0$; so $d$ is nilpotent iff there is some $N>0$ such that we can solve $a\sub{i}\sup{N}\equiv 0(n)$ for all $i\in\mb{N}$. \\[2mm]
 Given $\overline{a}\neq 0$ in $\mb{Z}\sub{n}$; let $1\leq a=bq\leq n-1$; where 
 $(q,n)=1$, i. e. $q$ is the coprime part of $a$, so
 $b=p\sub{1}\sup{m_{1}}\cdots p\sub{s}\sup{m_{s}}$ for some $m\sub{1},\dots,m\sub{2}$. 
  Then $\overline{a}$ is nilpotent in $\mb{Z}\sub{n}$ iff $p\sub{j}/a$ (or, equivalently, $1\leq m\sub{j}\leq r\sub{j}$)
 for each $j=1,\dots,s$. Since the multiplication on $\mb{Z}\sub{n}$ by a unit $\overline{q}$ is an isomorphism;
the order of $\overline{a}$ and $\overline{b}$ coincides, so we will assume without loss of generality that $q=1$.  \\[2mm]
 Assume that $d\neq0$ and suppose that the eigenvalues $\{\overline{a}\sub{i}: i\in\mb{N}\}$ of $d$
 satisfy those conditions. The representative elements $a\sub{i}$ 
 are taken from a finite subset of $\{1,\dots,n-1\}$. For each
 $i\in\mb{N}$ pick the first $k\sub{i}>0$ such that $a\sub{i}\sup{k_i}\equiv0(n)$. Since
 the set $\{k\sub{i}:i\in\mb{N}\}$ is finite; the order 
 of $d$ is $N=\max\{k\sub{i}:i\in\mb{N}\}$.\\[2mm]
  The generic model
  $\mc{Z}\sub{n}\mc{S}\sup{gen}=\sup{*}\!\!\!\mb{Z}\sub{n}$ is the set of non-standard $\mb{Z}\sub{n}$-sequences, i. e.
  the   ultrapower of $\mb{Z}\sub{n}$ modulo some ultrafilter $\mb{F}$ in $\mb{N}$. 
  We get a well defined   differential morphism $\mc{Z}\sub{n}\mc{S}\sup{gen}\ARROW{d}\mc{Z}\sub{n}\mc{S}\sup{gen}$
  given by $d([\overline{x}])=[d(\overline{x})]$ 
  where $[\overline{x}]$ is the germ of an element $\overline{x}\in\mc{Z}\sub{n}\mc{S}\sub{U}$, $U\subset\mb{N}$.
  Following \cite{dubois,kapranov} $\left(\sup{*}\!\mb{Z}\sub{n},d\right)$ is a $N$-complex. For $1\leq m\leq N-1$,
  the $m$th-amplitude generic cohomology is
  {\small\[
  	H\sub{m}\sup{gen}(X,\mc{Z}\sub{n}\mc{S}):=H\sub{m}\left(\sup{*}\!\mb{Z}\sub{n}\right)=
  	\frac{\text{ker}\left[\sup{*}\!\mb{Z}\sub{n}\ARROW{\hskip-2mm d^{k}}\ \sup{*}\!\mb{Z}\sub{n}\right]}{\text{Im}\left[\sup{*}\!\mb{Z}\sub{n}
  	\ARROW{\hskip-3mm d^{(N-k)}}\ \sup{*}\!\mb{Z}\sub{n}\right]}
  \]}
	Interesting examples arise as $n$ has many divisors. 
  For instance:
  \begin{itemize}
  		\item \underline{$n=12$}: Take $d(\overline{x})=(\overline{6x})$. Since 
  		$\overline{6}\sup{2}=\overline{36}=0$, $d\sup{2}=0$ and 
  		$d$ is a usual differential operator. 
  		$\ker(d)=\{\overline{x}\in \sup{*}\!\mb{Z}\sub{12} : x\sub{i}=\text{even }\forall i\}
  		\cong\ \sup{*}\!\mb{Z}\sub{6}$. Similarly, it can be seen that
  		$\text{Im}(d)\cong \sup{*}\!\mb{Z}\sub{2}$ so 
  		$H\sup{gen}(\mc{Z}\sub{12}\mc{S},d=\overline{6})=\frac{\sup{*}\!\mb{Z}\sub{6}}{\sup{*}\!\mb{Z}\sub{2}}=\ 
  		\sup{*}\!\mb{Z}\sub{3}$.
  		\item \underline{$n=48$}: The eigenvalues $\overline{a}\sub{i}$ of $d$ correspond to the integers
  		$\{\overline{a}\sub{i}: a\sub{i}=6,12,24\}$. 
  		The cohomology 	$H\sup{gen}(X,\mc{Z}\sub{48}\mc{S})$ will depend,  
		on the dimension of the subspaces of eigenvectors. Since at least one of them
		must have infinite  dimension, there always be a non-standard factor. 
  		Up to some adjusts, we can rearrange the canonical base as the disjoint union of\linebreak  
  		subspaces of eigenvectors. We obtain the following table of (usual) generic cohomologies ($N=2$):
  		{\small \begin{center}
	\begin{tabular}{|c|c|}\hline
  		$d$ &  $H\sup{gen}(X,\mc{Z}\sub{48}\mc{S})$ 
  		\\[1mm]  		
  		\hline
			$\overline{12}$  & $\sup{*}\!\mb{Z}\sub{3}$\\[1mm]
			$\overline{24}$  & $\sup{*}\!\mb{Z}\sub{12}$\\[1mm]
			$\overline{12}\oplus\overline{24}$  & $\mb{Z}\sub{3}\sup{\kappa}\oplus \sup{*}\!\mb{Z}\sub{12}$, or
			$\sup{*}\!\mb{Z}\sub{3}\oplus\mb{Z}\sub{12}\sup{\kappa}$, or $\sup{*}\!\mb{Z}\sub{3}\oplus\sup{*}\!\mb{Z}\sub{12}$ \\[1mm]
  		\hline
	\end{tabular}\vskip2mm
	\end{center}}
	In the last row, the factor with a superscript $\kappa\in\mb{N}$ corresponds to a finite-dimensional
	subspace of eigenvectors, whenever it is the case.	
	The order of the (diagonal) differential operator $d$ induced by the eigenvalue $\overline{6}$ is $N=4$. 
	The amplitude generic cohomology $H\sub{m}\sup{gen}(X,\mc{Z}\sub{48}\mc{S})$
	for $m=1,2,3$ is always $\sup{*}\!\mb{Z}\sub{3}$. Finally; the differential operators
	$\overline{6}\oplus\overline{12}$, $\overline{6}\oplus\overline{24}$ and
	$\overline{6}\oplus\overline{12}\oplus\overline{24}$ can be treated in a similar way.
  \end{itemize}

\subsection{A final comment}\label{ssection a final comment} 
A description of $H\sup{gen}$ as a loose  cohomology and its relation to Weil cohomologies,  
in the terms of Macintyre \cite{macintyre2}; is still a pending assingment.
The above examples are similar to those non standard cohomologies of Br\"unjes and Serp\'e, 
\cite{brunjes}. New examples of generic cohomologies should appear from
structure sheaves such as those provided by Abramsky \cite{abramsky}. 
An additional reason of interest in this subject, from a model
theoretic perspective, is the development of stability theoretical
tools for the classification of sheaves. These question are not treated
here. The authors hope to fulfill these tasks in a forthcoming article.

\end{document}